\documentclass[10pt]{article}
\usepackage{amssymb,amsthm,amsmath}
\usepackage{color}
\theoremstyle{definition}
\newtheorem{definition}{Definition}
\newtheorem{example}{Example}
\newtheorem{remark}{Remark}
\theoremstyle{theorem} 
\newtheorem{theorem}{Theorem}
\newtheorem{proposition}{Proposition}
\newtheorem{lemma}{Lemma}
\newtheorem{corollary}{Corollary}

\title{Multivariable codes over finite chain rings: semisimple codes}
\author{E. MART\'INEZ-MORO\thanks{Departamento de Matem\'atica Aplicada, Universidad de Valladolid, Spain, ({\tt edgar@maf.uva.es}). Partially supported by MEC MTM2004-00876 and MTM2004-00958 I+D projects.} \and I.F. R\'UA\thanks{Departamento de Matem\'aticas, Estad\'{\i}stica y Computaci\'on
Universidad de Cantabria, Spain, ({\tt i.f.rua@unican.es}).
Partially supported by MTM2004-08115-C04-01 I+D project.}}
\begin{document}

\maketitle

\begin{abstract}
The structure of multivariate semisimple codes over a finite chain
ring $R$ is established using the structure of the
residue field $\bar R$. Multivariate codes extend in a natural way
the univariate cyclic and negacyclic codes and include some
non-trivial codes over $R$. The structure of the dual codes in the
semisimple abelian case is also derived and some conditions on the
existence of selfdual codes over $R$ are studied.\\[0.5em]



\noindent Submitted to:  \textit{SIAM Journal on Discrete Mathematics}.
\end{abstract}

\section{Introduction}
Many authors have stated that  many classical codes are ideals in certain algebras over a finite field, see for example \cite{Berman,Charpin,Poli2}. On the other hand, the theory of error-correcting codes over finite rings has gained a great relevance since the realization that some non-linear codes can be seen as linear codes over a finite ring (see for example  \cite{Calderbank1,Calderbank2,Nec82,Nechaev,NecKuz96}). This paper is a contribution to both lines pointed above and its purpose is to describe multivariate semisimple codes over a finite chain ring $R$. Through the paper a semisimple code over $R$ will be an ideal of a particular type of $R$-algebras. We shall note that the name of semisimple codes arise from the fact that the image code in the residue ring $\bar R$ is semisimple (in fact they are not semisimple over $R$). The main tools used in the paper are Hensel's Lemma and the decomposition of the roots of the defining ideal in cyclotomic classes. Multivariate codes extend in a natural way the univariate cyclic and negacyclic codes \cite{CN-FCR}  and include some non-trivial codes over $R$.

The paper is organized as follows. In Section 2 we present the basic results on finite chain rings needed. Section 3 is devoted to the definition of the codes and their ambient space as well as the description of their structure. In Section 4  we study the duals of abelian semisimple codes. Finally in Section 5 we characterize those non-trivial abelian semisimple codes that are self-dual.

\section{Preliminaries}

In this section we fix our notation and show some basic facts about
finite chain rings (see for example \cite{B-F,McD} for a complete
account). From now on, by a ring $R$ we will always mean an
associative commutative ring with identity, unless explicitly
stated.  A ring  $R$ is called \textbf{local ring}
if it has a unique maximal ideal and
 it is called a \textbf{chain ring}  if the set of all the ideals is a chain under set-theoretical inclusion.
 It can be shown (see for example Proposition 2.1 in \cite{CN-FCR}) that  $R$ is a finite commutative chain ring if,
  and only if, $R$ is a local ring and its maximal ideal $M$ is principal. In this case let $a\in R$ be a fixed generator of the ideal $M=\mathrm{rad}(R)$ and, since $a\in M$ is nilpotent, let $t$ be its nilpotency index. Then we have
\begin{equation}
\left\langle 0\right\rangle =\left\langle a^t \right\rangle \subsetneq \left\langle a^{t-1} \right\rangle \subsetneq \dots \subsetneq \left\langle a^1 \right\rangle= M \subsetneq \left\langle a^0 \right\rangle=R.
\end{equation}
Let $q=p^l$ where $p$ is a prime and  $\mathbb F_q=\bar R=R/M$ is the residue field of $R$. We can extend the natural ring homomorphism $r\mapsto \bar r= r+M$ as follows
\begin{equation}
\begin{array}{rcl}
R & \hookrightarrow &  R[X]\\
\bar{} \downarrow & & \bar{} \downarrow \\
\mathbb F_q & \hookrightarrow  & \mathbb F_q[X]
\end{array}
\end{equation}
Two polynomials $f_1,f_2\in R[X]$ are \textbf{coprime} if $(f_1,f_2)=1$.  A polynomial $f\in R[X]$ is called \textbf{regular} if it is not a zero divisor and \textbf{basic irreducible} if it is regular and $\bar f\in \mathbb F_q[X]$ is irreducible. The following well known result will be used several times in the paper, for a proof see for example \cite[Theorem 3.2.6]{B-F}

\begin{theorem}[Hensel's lemma] Let $R$ be a finite local ring and $f\in R[X]$ be a monic polynomial
such that $\bar f=g_1g_2\dots g_r$ where the
polynomials $g_i\in \bar R[X]$ are monic and
pairwise relative prime. Then there exist monic coprime polynomials
$f_i\in R[X]$  $i=1,\ldots, r$ such that $\bar f_i=gi$ for all
$i=1,\ldots, r$ and $f=f_1f_2\dots f_r$. This decomposition is
uniquely determined up to a permutation of the factors.
\end{theorem}

From Hensel's lemma we can deduce the existence of polynomials \textbf{ lifting}  a factorization in $\bar R[X]$ to a factorization in $R[X]$.  We refer to these polynomials as \textbf{lifting factors}.

Let $R$ and $S$ be two rings such that $R\subseteq S$, then we say that $S$ is an \textbf{extension} of $R$. If $T\subseteq S$ and $T\neq\emptyset$ of finite cardinality, then the ring \textbf{generated} by $T$ is the smallest subring $A$ containing $R \cup T$. If $T=\{a\}$ is a singleton, then we call the extension \textbf{simple} and denote it by $A=R(a)$. If $R$ and $S$ are two finite local rings with residue fields $F$ and $K$ respectively, such that $R\subseteq S$, then $S$ is a \textbf{separable extension} of $R$ if $K$ is a separable extension of $F$ in the sense of field extensions.

In our paper we consider monic polynomials
$t_i(X_i)\in R[X_i]$ $i=1,\ldots ,r$ such that
$\bar{t_i}(X_i)\in K[X_i]$ is
 square-free, where $K$ is the algebraic closure of $\mathbb F_q$ (semisimple case). So
we have that
\begin{equation}
t_i(X_i)=
\prod_{j=1}^{r_i}f_{i,j}(X_i)
\end{equation}
where $f_{i,j}(X_i){,}\;\ j=1,\dots,r_i$ are monic basic irreducible
polynomials and  $(f_{i,j},f_{i,k})=1$, if $j\not=k$.
 This decomposition is
unique up to a relabelling of the factors {due to Hensel's lemma}.

\section{Multivariable semisimple codes}

In this section we will obtain the structure of a multivariable
semisimple code over a finite chain ring $R$, {i.e.}, we will
describe explicitly the structure of the ideals of the ring
$R[X_1,\dots,X_r]/\left\langle t_1(X_1),\dots,t_r(X_r)\right\rangle
$. In order to obtain this description we will decompose this ring
as a direct sum of finite local chain rings. This decomposition is
based on the corresponding decomposition of the semisimple ring
$\mathbb F_q[X_1,\dots,X_r]/\left\langle\overline
t_1(X_1),\dots,\overline t_r(X_r)\right \rangle$.

\subsection{Descomposition of $R[X_1,\dots,X_r]/\left<  t_1(X_1),\ldots ,t_r(X_r)\right>$}

Let $$I=\left<  t_1(X_1),\ldots ,t_r(X_r)\right>\lhd R[X_1,\ldots ,X_r]$$ be the ideal generated by the polynomials $t_i(X_i)$ $i=1,\ldots ,r$ defined as in the section above. Let $H_i$ be  the set of roots of $\bar t_i(X_i)$ in an suitable extension field of $\mathbb F_q$ for each $i=1,\ldots ,r$ (notice that
$\bar t_i(X_i)$ has no multiple roots).
\begin{definition} Let $\mu=(\mu_1,\ldots ,\mu_r)\in H_1\times \ldots \times H_r$,
then we define the \textbf{class} of  $\mu$ as
\begin{equation}
C(\mu)=\left\lbrace (\mu_1^{q^s},\ldots ,\mu_r^{q^s})\mid s\in\mathbb N\right\rbrace .
\end{equation}
\end{definition}

\begin{proposition}
Let $\mu=(\mu_1,\ldots ,\mu_r)\in H_1\times \ldots \times H_r$ and
$d_i$ be the degree of the minimal polynomial of $\mu_i$
over $\bar R=\mathbb F_q$ for each $i=1,\ldots
,r$, then we have that
\begin{enumerate}
\item $\vert C(\mu)\vert =\mathrm{l.c.m.}(d_1,d_2,\ldots ,d_r)=[\mathbb F_q(\mu_1,\ldots ,\mu_r):\mathbb F_q]$.
\item The set of classes $C(\mu)$ is a partition of $H_1\times \ldots \times H_r$.
\item For each ideal $J\lhd \mathbb F_q[X_1,\ldots ,X_r]/\left<  \bar t_1(X_1),\ldots ,\bar t_r(X_r)\right>$ the affine variety $V(J)$ of common zeros of the elements in $J$  is a union of classes.
\end{enumerate}
\end{proposition}

\begin{proof}
See \cite{Poli} for a proof.
\end{proof}

\begin{definition}\label{def:pw}
Let us denote by $\mathrm{Irr}(\alpha ,\mathbb F_q)$ the minimal polynomial of $\alpha\in K$ over the field $\mathbb F_q$ ($K$ is an algebraic extension of $\mathbb F_q$). If $\mu=(\mu_1,\ldots ,\mu_r)\in H_1\times \ldots \times H_r$, then we define the following polynomials:
\begin{enumerate}
\item $p_{\mu,i}(X_i)=\mathrm{Irr}(\mu_i ,\mathbb F_q)$, and $d_{\mu,i}=\deg p_{\mu,i}$ for all $i=1,\ldots ,r$.
\item $w_{\mu,i}(\mu_1,\ldots ,\mu_{i-1}, X_i)= \mathrm{Irr}(\mu_i ,\mathbb F_q(\mu_1,\ldots ,\mu_{i-1}))$ for all $i=2,\ldots ,r$.
\item $\pi_{\mu,i}(\mu_1,\ldots ,\mu_{i-1}, X_i)=p_{\mu,i}(X_i)/w_{\mu,i}(\mu_1,\ldots ,\mu_{i-1}, X_i) $ for all $i=2,\ldots ,r$.
\end{enumerate}
\end{definition}

\begin{remark}\label{rem:1}
All the polynomials in the definition above can be seen as polynomials in $\mathbb F_q[X_1,\ldots ,X_r]$ (substituting $\mu_i$ by $X_i$) and clearly the following ring isomorphism holds
\begin{equation}
\mathbb F_q[X_1,\ldots ,X_r]/\left\langle
p_{\mu,1},w_{\mu,2},\ldots ,w_{\mu,r}\right\rangle\cong \mathbb
F_q(\mu_1,\ldots ,\mu_r).
\end{equation}
Moreover, if $\mu'\in C(\mu)$, then $p_{\mu,i}=p_{\mu',i}\ i=1,\dots,r$ and $w_{\mu,i}=w_{\mu',i},\pi_{\mu,i}=\pi_{\mu',i}\ i=2,\dots,r$.
\end{remark}

If $q(X)\in R[X]$ is the Hensel's lifting of a monic irreducible polynomial $p(X)\in\mathbb F_q[X]$ and $M=\hbox{rad}(R)$, then $\left\langle M,q(X)\right\rangle $ is a maximal ideal of $R[X]$ and (cf. \cite[Remark after Lemma 3.2.10]{B-F})

\[
R[X]/\left\langle M,q(X)\right\rangle  \cong \mathbb F_q[X]/\left\langle p(X)\right\rangle \cong \mathbb F_q(\alpha)
\]
where $p(\alpha)=0$. Notice that $S=R[X]/\left
< q(X)\right>$ is a local ring with maximal ideal $\left
<M,q(X)\right
>+\left<q(X)\right >$, that can be seen as a separable extension of $R$ (since $p(X)\in\mathbb F_q[X]$ is
irreducible). In particular we have that $S$ is a finite local chain
ring. If we consider $q(X)\in S[X]$, then the element
$A=X+\left<q(X)\right >\in S$ is a root of the polynomial $q(X)$
that lifts $\alpha$, and so we can write $S=R(A)$.

\begin{definition}\label{def:qz} Let $\mu$, $R$, $p_{\mu,i}\ i=1,\dots,r$, $w_{\mu,i}$ and $\pi_{\mu,i}$ $i=2,\ldots ,r$ be as in Definition \ref{def:pw}, then for all $i=1,\dots,r $ we define $q_{\mu,i}$ as the Hensel's lifting of the polynomial $p_{\mu,i}$ to $R[X_i]$ and, for all $i=2,\dots,r$, we define $z_{\mu,i}$ and $\sigma_{\mu,i}$ as the Hensel's liftings of the polynomials $w_{\mu,i},\pi_{\mu,i}\in \mathbb F_q(\mu_1,\ldots,\mu_{i-1})[X_i]$ to  $R_{i-1}[X_i]$ where $R_{i-1}$ is the local ring $R(\mu_1, \ldots ,\mu_{i-1})$.
\end{definition}

\begin{remark} By the discussion above the polynomials $z_{\mu,i}$ and $\sigma_{\mu,i}$ $i=2,\ldots ,r$ are well defined.
  Moreover, as in Remark \ref{rem:1} they can be seen as polynomials in $ R[X_1,\ldots ,X_r]$ (substituting the lifting
 of the root $\mu_i$ by the corresponding indeterminate $X_i$), and
 $T= R[X_1,\ldots ,X_r]/\left\langle q_{\mu,1},z_{\mu,2},\ldots ,z_{\mu,r}\right\rangle $
 is a local ring with maximal ideal
 $\mathfrak{m}=\left\langle M,q_{\mu,1},z_{\mu,2},\ldots ,z_{\mu,r}\right\rangle +\left\langle q_{\mu,1},z_{\mu,2},\ldots ,z_{\mu,r}\right\rangle $
 and quotient ring
\begin{equation}
T/\mathfrak{m}\cong \mathbb F_q(\mu_1,\ldots
,\mu_r).
\end{equation}
\end{remark}
\begin{lemma}\label{lem:1}
Let $R$ be a finite chain ring with maximal ideal $M=\left\langle a
\right\rangle$ and residue field $\mathbb F_q$ where  the nilpotency
index of $a$ is $t$.  Let $\mu=(\mu_1,\ldots ,\mu_r)\in H_1\times
\ldots \times H_r$ and consider the ideal
\begin{equation}
I_{\mu}=\left\langle q_{\mu,1},z_{\mu,2},\ldots ,z_{\mu,r}\right\rangle
\end{equation}
where the polynomials $q_{\mu,1},z_{\mu,i}$ $i=2,\ldots , r$ are
defined as above.

\noindent Then $R[X_1,\ldots , X_r]/I_{\mu}$ is a finite commutative
chain ring with maximal ideal $\left\langle a+I_{\mu}\right\rangle
$, residue field $\mathbb F_q (\mu_1,\ldots ,\mu_r)$ and precisely
the following ideals
\begin{equation}
\left\langle 0\right\rangle =\left\langle a^t+ I_{\mu}\right\rangle \subsetneq \left\langle a^{t-1}+ I_{\mu} \right\rangle \subsetneq \dots \subsetneq \left\langle a^1 + I_{\mu} \right\rangle= M \subsetneq \left\langle a^0 + I_{\mu}\right\rangle.
\end{equation}
\end{lemma}
\begin{proof}
It is a straightforward  conclusion of the above discussion and the fact that  $M=\left\langle a \right\rangle$
\end{proof}

\begin{definition}
 Let $\mu=(\mu_1,\ldots ,\mu_r)\in H_1\times \ldots \times H_r$, we
define the following polynomial in $R[X_1,\ldots ,X_r]$
 \begin{equation}  h_\mu(X_1,\ldots ,X_r)=\prod_{i=1}^{r}
\frac{t_i(X_i)}{q_{\mu,i}(X_i)}\prod_{i=2}^r
\sigma_{\mu,i}(X_2,\ldots ,X_r)  \end{equation}  where the
polynomials $t_i,q_{\mu,i}$ $i=1,\ldots , r$ and $\sigma_{\mu,i}$
$i=2,\ldots , r$ are defined as in {Definition \ref{def:qz}}.
\end{definition}

\begin{proposition}\label{prop:ann} If $I=\left<  t_1(X_1),\ldots ,t_r(X_r)\right> \lhd R[X_1,\ldots ,X_r]$, then the
annihilator of $\left\langle h_\mu+I\right\rangle $  in $R[X_1,\ldots ,X_r]/I$ is
\begin{equation}
\mathrm{Ann}\left( \left<h_\mu+I\right>\right) =I_\mu +I
\end{equation}
\end{proposition}
\begin{proof}
Clearly $I_\mu +I \subseteq \mathrm{Ann}\left( \left<h_\mu+I\right>\right)$.

 \noindent On the other hand,
if $g+I\in \mathrm{Ann}\left(\left<h_\mu+I\right>\right)$, then the
polynomial $\overline{gh}_\mu\in \bar I=\left\langle \bar
t_1(X_1),\dots,\bar t_r(X_r)\right\rangle $ and so $\bar g+\bar I\in
\mathrm{Ann}\left(\left< \bar h_\mu+\bar I\right >\right)=\left<\bar
q_{\mu,1},\bar z_{\mu,2},\dots,\bar z_{\mu,r}\right>$ (cf.
\cite[Proposition 6]{Poli}). Hence $g+I\in
\left\langle I_\mu+\left\langle a\right\rangle \right\rangle+I$ and
{thus}
$\mathrm{Ann}\left\langle\left<h_\mu+I\right>\right\rangle=\left\langle
I_\mu+\left\langle a^s\right\rangle \right\rangle+I$ for some
$s\in\{0,\dots,t\}$. Now, if ${\theta_i}$ is a root of $q_{\mu,i}\
i=1,\dots,r$ lifting $\mu_i$ and we denote $\Theta=({
\theta_1,\dots,\theta_r})$, then $h_\mu(\Theta)\not \in
\left<a\right >$ (since $\bar h_\mu(\mu)\not =0$, cf. \cite[Chapter
5, Proposition 7]{Poli}) and {therefore} we can conclude
$\mathrm{Ann}\left(\left<h_\mu+I\right>\right)=I_\mu+I$ as desired
(otherwise $s<t$, and so $a^{t-1}=a^sa^{t-1-s}\in
\mathrm{Ann}\left(\left<h_\mu+I\right>\right)$ implies
$a^{t-1}h_\mu\in I$ and $0=a^{t-1}h_\mu(\Theta)$, i.e.,
$h_\mu(\Theta)\in \left\langle a\right\rangle $, a contradiction).
\end{proof}

Notice that, if $\mu'\in C(\mu)$, then $q_{\mu,i}=q_{\mu',i}\ i=1,\dots,r$, $z_{\mu,i}=z_{\mu',i},\sigma_{\mu,i}=\sigma_{\mu',i}\ i=2,\dots,r$ and so $h_{\mu}=h_{\mu'}$. Therefore, by abuse of notation we shall write $I_C$ and $h_C$ instead of $I_\mu$ and $h_\mu$ provided that $C$ is the class $C(\mu)$.

\begin{lemma}\label{lemma:2} Let $\mathcal C$ be the set of classes $C(\mu)$ where $\mu\in H_1\times \ldots \times H_r${, and } ${ C,C^\prime\in \mathcal C}$.  Then:
\begin{enumerate}
\item The set of zeros of $\bar h_{C}$ is $ H_1\times \ldots \times H_r\setminus C$
 and the set of zeros of $\bar I_C$ is $C$.
\item $\left\langle t_1(X_1),\dots,t_r(X_r)\right\rangle =\bigcap_{C\in \mathcal C} I_C$.
\item $I_C$, $I_{C^\prime}$ are comaximal if $C\neq C^\prime$.
\end{enumerate}
\end{lemma}

\begin{proof}$ $
\begin{enumerate}
\item Is a direct translation of Proposition 7 in \cite[Chapter 5]{Poli}. Note that the ideal $\bar I=\left\langle \bar t_1(X_1),\dots,\bar t_r(X_r)\right\rangle$ is a radical ideal in $\bar {\mathbb F}_q[X_1,\ldots, X_r]$ and the variety
\begin{equation}\label{eq:disun} V(\left\langle \bar t_1(X_1),\dots,\bar t_r(X_r)\right\rangle)=\bigsqcup_{C\in \mathcal C} C= V(\bar I_C)
\end{equation}
thus $\left\langle \bar t_1(X_1),\dots,\bar t_r(X_r) \right\rangle =\bigcap_{C\in \mathcal C} \bar I_C$.
\item {Clearly $\left\langle t_1(X_1),\dots,t_r(X_r)\right\rangle \subseteq\bigcap_{C\in \mathcal C} I_C$. Suppose that $f\in\bigcap_{C\in \mathcal C} I_C$, then by Proposition \ref{prop:ann} we have that $f+I\in \mathrm{Ann}\left( \left<h_\mu+I\right>\right)$ for all choices of $\mu$. Thus $\overline{fh_\mu}\in \bar I$ for all $\mu$, and by part 1) of this proof $\bar f\in \bar I$ and the result follows. }

\item Arises from the fact that in equation (\ref{eq:disun}) the union is disjoint.
\end{enumerate}
\end{proof}

\begin{theorem}\label{thm:CRT}
\begin{equation} R[X_1,\dots,X_r]/I\cong \bigoplus_{C\in \mathcal C} \left\langle h_C+I\right\rangle\end{equation}
where $\left\langle h_C+I\right>\cong R[X_1,\dots,X_r]/I_C$ is a finite commutative chain ring whith maximal ideal $\left\langle a+I_C\right\rangle$.
\end{theorem}
\begin{proof}
By the Chinese Remainder theorem
\[
R[X_1,\dots,X_r]/I=R[X_1,\dots,X_r]/\bigcap_{C\in \mathcal C} I_C\cong\bigoplus_{C\in \mathcal C} R[X_1,\dots,X_r]/I_C
\]
and the result follows.
\end{proof}

\begin{remark}\label{idempotentes}
The above theorem is equivalent to the fact that there exist
primitive orthogonal idempotents elements $e_i\in
R[X_1,\dots,X_r]/I$ (one for each class $C_i\in \mathcal C$) such
that  $1=\sum e_i$ and $e_i\left( R[X_1,\dots,X_r]/I \right)\cong
\left\langle h_{C_i}+I\right\rangle$ (cf. \cite[Proposition
3.1.3]{B-F}). Namely, the idempotent $e_i$ is exactly the element
$g_{C_i}h_{C_i}+I$, where
$g_{C_i}h_{C_i}+I_{C_i}=1+I_{C_i}$.
\end{remark}

\subsection{Description of the codes}

Classical coding theory has been developed in vector spaces over
finite fields, a good background in algebraic codes over finite
fields is the textbook \cite{M-S}.  We describe some natural
modifications that leads us to codes over finite rings, see for
example the textbook \cite{B-F}.

For a finite {commutative ring $R$ consider the set $R^n$ of all
$n$-uples as a module over $R$ as usual. We say that a subset
$\mathcal K$ of $R^n$ is a \textbf{linear code} if $\mathcal K$ is
an $R$-submodule of $R^n$. Given an ideal $J\lhd R[X_1,\dots,X_r]$
such that the algebra $R[X_1,\dots,X_r]/J$ has
finite rank $n$ as $R$-module, and given an
ordering on the set of terms, each element of $R[X_1,\dots,X_r]/J$
can be identified with a $n$-uple in $R^n$.}

 Given two elements $\mathbf{x}=(x_1,\ldots ,x_n),\mathbf{y}=(y_1,\ldots ,y_n)\in R^n$ the scalar product is   $\mathbf{x}\cdot \mathbf{y}=(x_1y_1+\ldots +x_ny_n)\in R$. We say
{that $\mathbf{x},\mathbf{y}$ are \textbf{orthogonal} if $\mathbf{x}\cdot\mathbf{y}=0$ and, for a linear
  code $\mathcal K$, we define the \textbf{dual code} as $\mathcal K^\perp=\left\lbrace \mathbf{x}\in R^n\mid
  \mathbf{x}\cdot\mathbf{c}=0\quad \forall \mathbf{c}\in \mathcal K \right\rbrace$. The code $\mathcal K$
  is called
  \textbf{selfdual} if $\mathcal K=\mathcal K^\perp$.}

 \begin{definition}[{\bf Multivariable semisimple code}]\label{def:sem-code} Let $t_i(X_i)\in R[X_i]$ $i=1,\ldots ,r$ be polynomials
 over a finite chain ring $R$. A \textbf{multivariable code} is an ideal $\mathcal K$ of the ring $R[X_1,\dots,X_r]/\left\langle t_1(X_1),\dots,t_r(X_r)\right\rangle $. If the polynomials $t_i$, $i=1,\dots,r$, are defined as in the previous section, then we shall say that the code is \textbf{semisimple}.
\end{definition}

Notice that a multivariable semisimple code is not \emph{semisimple} in the classical ring theoretic sense. Indeed, we shall see later (Corollary \ref{cor:des-ideales}) that any semisimple code is a sum of finite chain rings. The name is justified so, by the fact that the image code $\overline {\mathcal K}$ of $\mathcal K$ in $\overline R[X_1,\dots,X_r]/\left\langle \overline t_1(X_1),\dots,\overline t_r(X_r)\right\rangle$ is semisimple ($\overline{\mathcal K}$ is a sum of simple ideals).

Clearly this class of codes includes, among others,
cyclic and negacyclic semisimple codes.
Next we present an example of non-trivial codes that fall into this category. This example is due to A.A. Nechaev and
A.S. Kuzmin \cite{Nechaev97}.

\begin{example}
    Let $R=GR(q^2,2^2)$ ($q=2^l$) be the Galois Ring of cardinality $q^2$ and characteristic $2^2$ \cite{McD}, and let
    $S=GR(q^{2m},2^2)$ be its Galois extension of odd degree $m\ge 3$. Both $R$ and $S$ are finite commutative chain rings
     with maximal ideals $2R$ and $2S$ and residue fields $\overline{R}=GF(q)$ and
     $\overline{S}=GF(q^m)$, respectively. With the help of the Teichm\"uller Coordinate Set (TCS)
     $\Gamma(S)=\{a^{q^m}=a\ |\ a\in S\}$ any element $a\in S$ can be decomposed uniquely
     as $a=\gamma_0(a)+2\gamma_1(a)$, where $\gamma_i(a)\in \Gamma(S)$. Moreover, if $\oplus:\Gamma(S)\times \Gamma(S)\to \Gamma(S)$ is defined as $a\oplus b=\gamma_0(a+b)$, then $(\Gamma(S),\oplus,\cdot)$ is the finite field $GF(q^m)$ whose cyclic multiplicative group is generated by an element $\theta$ of order $\tau=q^m-1$, and the TCS $\Gamma(R)=\{a^q=a\ |\ a\in R\}=\{w_0=0,w_1,\dots,w_{q-1}\}$ is the subfield $GF(q)$. Let $\hbox{Tr}:S\to R$ denote the trace function from $S$ onto $R$, then the (shortened) \emph{$R$-base linear code}  is given by:
$$\mathcal L=\{(\hbox{Tr}(\xi)+a,\hbox{Tr}(\xi\theta)+a,\dots,\hbox{Tr}(\xi\theta^{\tau-1})+a)\ |\ \xi\in S, a\in R\}.$$
It is an $R$-linear code of length $\tau$, cardinality $q^{2(m+1)}$ and
the (shortened) \emph{Generalized Kerdock code} is the projection of $\mathcal L$ in $\Gamma(R)^{\tau q}$ with the help of $\tau$ copies of the $RS$-map:
$$\gamma_*(a)=(\gamma_1(a),\gamma_1(a)\oplus w_1\gamma_0(a),\dots,\gamma_1(a)\oplus w_{q-1}\gamma_0(a))\ ,\ a\in R.$$
It is an $GF(q)$-nonlinear code of length $\tau q$, cardinality $q^{2(m+1)}$ and Hamming distance $\frac{q-1}{q}(n-\sqrt n)-q$.

This code can be presented in a polycyclic form with the help of a
multivariable code over the finite chain ring $R$, by the following
way. The multiplicative group $U=1+2R=\{u_0=1,u_1,\dots,u_{q-1}\}$
is a direct product $<\eta_1>\times\dots\times<\eta_l>$ of $l$
subgroups of order $2$. Consider the ideal $I$ of
$R[X_1,\dots,X_r]$, where $r=l+1$, generated by the polynomials
$t_1(X_1)=X_1^{{\tau}}-1,t_2(X_2)=X_2^2-1,\dots,t_r(X_r)=X_r^2-1$.
If we denote $\overrightarrow U=(u_0,\dots,u_{q-1})$ and
$\overrightarrow a\otimes \overrightarrow U=(a_1\overrightarrow
U,\dots,a_q\overrightarrow U)\in R^{q\tau}$ for any $\overrightarrow
a\in R^\tau$, then the multivariable code $\mathcal K\lhd
R[X_1,\dots,X_r]/I$ given by
$$\mathcal K=\left\{\sum_{i_1=0}^{\tau-1}\sum_{i_2=0}^1\dots\sum_{i_r=0}^1\left((\hbox{Tr}(\xi\theta^{i_1})+a)\eta_1^{i_2}\dots\eta_l^{i_r}\right)
X_1^{i_1}X_2^{i_2}\dots X_r^{i_r}\ |\ \xi\in S, a\in R\right\}$$
is equivalent to the code $\mathcal L\otimes \overrightarrow U$, and the shortened Generalized Kerdock code is equivalent to the polycyclic code $\gamma_1^{q\tau}(\mathcal K)$. {Notice that this code is not semisimple, though.}
\end{example}

{Now we can back to the description of multivariable semisimple codes. The following two results
are straight forward corollaries of Theorem \ref{thm:CRT}.}

\begin{corollary}\label{cor:des-ideales} Let $R$ be a finite chain ring with maximal ideal $\left\langle a \right\rangle$
and nilpotency index $t$. Any semisimple code $\mathcal K$  in
$R[X_1,\dots,X_r]/I$ where $I=\left\langle
t_1(X_1),\dots,t_r(X_r)\right\rangle $, is a  sum of ideals of the
form
\begin{equation}
\left\langle a^{j_C}h_C+I\right\rangle \qquad 0\le j_C\le t,\hbox{
and } C\in\mathcal C
\end{equation}
\end{corollary}

\begin{corollary} {In the conditions of the previous corollary, there are $(t+1)^N$ semisimple codes in
$R[X_1,\dots,X_r]/I$, where $N=|\mathcal C|$.}
\end{corollary}

We shall now obtain an explicit description of
semisimple codes in terms of polynomials of the ring
$R[X_1,\dots,X_r]$.

\begin{theorem}\label{H}
{If $\mathcal K$ is a semisimple code in $R[X_1,\dots,X_r]/I$, then
there exists a family of
 polynomials  $G_0,\dots,G_t\in R[X_1,\dots,X_r]$ determining uniquely the ideals $\left\langle G_i+I\right\rangle$
 such that
 \begin{equation}
 I=\bigcap_{i=0}^t\mathrm{Ann}\left\langle G_i+I\right\rangle,\quad \mathcal K=\left\langle  G_1,a G_2,\dots,a^{t-1}
 G_t\right\rangle+I
 \end{equation}
 and, for each pair $0\leq i< j\leq t$, the ideals  $\mathrm{Ann}\left\langle G_i+I\right\rangle$, $\mathrm{Ann}\left\langle
 G_j+I\right\rangle$ are comaximal. Moreover, $\mathcal K=\left\langle
 G+I\right\rangle$, where $G=\sum_{i=0}^{t-1}a^iG_{i+1}$.}
\end{theorem}
\begin{proof}
By Corollary \ref{cor:des-ideales} $\mathcal K$ is
a direct sum of ideals of the form $\left\langle a^{j_C}h_C+I\right
\rangle$, where $0\le j_C\le t$, and $C\in\mathcal C$. If
$N=|\mathcal C| $ is the number of classes in $\mathcal C$, then,
after reordering of the classes in $\mathcal C$, we have
\begin{eqnarray*}
\mathcal K & = &   \left\langle h_{C_{k_1+1}}+I\right\rangle \oplus \dots
\oplus \left\langle h_{C_{{k_1}+{k_2}}}+I\right\rangle \\  & & \oplus
\left\langle ah_{C_{{k_1}+{k_2}+1}}+I\right\rangle \oplus \dots
\oplus
 \left\langle ah_{C_{{k_1}+{k_2}+{k_3}}}+I\right\rangle\oplus
\dots \oplus \\ & & \left\langle
a^{t-1}h_{C_{\sum_{i=1}^t{k_i}+1}}+I\right\rangle \oplus \dots\oplus
\left\langle a^{t-1}h_{C_{N}}+I\right\rangle
\end{eqnarray*}

where $k_i\geq 0$ for all $i=1,2,\ldots ,t$ and
$\sum_{i=1}^t{k_i}+1\leq N$. Let $k_0=0$ and
$k_{t+1}=N-\sum_{i=1}^t{k_i}$, and define $$G_i=
\sum_{j=k_0+\dots+k_i+1}^{k_0+\dots+k_{i+1}}g_{C_j}h_{C_j}$$ where
$g_{C_j}\in R[X_1,\dots,X_r],\
j=k_0+\dots+k_i+1,\dots,k_0+\dots+k_{i+1},\ i=0,\dots,t$ are the
polynomials defining the primitive orthogonal idempotents of Remark
\ref{idempotentes}. Then:
$$\left\langle
G_i+I\right\rangle=\sum_{j=k_0+\dots+k_i+1}^{k_0+\dots+k_{i+1}}\left\langle
h_{C_j}+I\right \rangle$$
and so we have $\mathcal K=\left\langle  G_1,a
G_2,\dots,a^{t-1}
 G_t\right\rangle+I$, and
$$\bigcap_{i=0}^t\mathrm{Ann}\left\langle G_i+I\right\rangle
={\bigcap_{i=0}^t\bigcap_{{j=k_0+\dots +k_i+1}}^{k_0+\dots
+k_{i+1}}}\mathrm{Ann}\left(\left<h_{C_j}+I\right>\right)=\bigcap_{{k=0}}^{N}I_{C_k}+I.$$
Moreover, for each pair $0\leq i< j\leq t$, the
ideals $\mathrm{Ann}\left\langle G_i+I\right\rangle$,
$\mathrm{Ann}\left\langle
 G_j+I\right\rangle$ are comaximal,
 from \textit{2)} and \textit{3)} in Lemma
\ref{lemma:2}. The uniqueness of the ideals $\left\langle
G_i+I\right\rangle, i=0,\dots,t$, follows from fact that the
decomposition in Theorem \ref{thm:CRT} is unique, and Corollary
\ref{cor:des-ideales}. {Finally, the equality $\mathcal
K=\left\langle G+I\right\rangle$ is satisfied, since each elements
$G_i$ is a sum of primitive idempotent orthogonals of the ring.}
\end{proof}

With this description in hand we can obtain the cardinality of any
semisimple code.

\begin{corollary}
In the conditions of Theorem \ref{H} $R[X_1,\dots,X_r]/I$ is a
principal ideal ring and, for any semisimple code $\mathcal K$, we
have:
$$|\mathcal K|=|\bar R|^{\sum_{i=0}^{t-1}(t-i)N_i} $$
where $N_i$ denotes the number of zeros $\mu\in H_1\times\dots H_r$
of $\bar G_i, i=0,\dots, t-1$.
\end{corollary}

\begin{proof}
For $i=0,\dots,t-1$ we have
$$\left\langle a^iG_{i+1}+I\right\rangle=\left(\frac{|R|}{|\left\langle a^{i}\right\rangle|}\right)^{\hbox{rank}_R
(\left\langle G_i+I\right\rangle)} =|\bar R|^{(t-i)\hbox{rank}_R
(\left\langle G_i+I\right\rangle)}.$$

Since $\hbox{rank}_R(\left\langle
G_i+I\right\rangle)=\hbox{dim}_{\bar R}\left\langle \bar G_i+\bar
I\right\rangle$, the result follows from \cite{Poli}.
\end{proof}

\subsection{Hamming distance of the codes}

For $\mathbf{c}\in R^n$ we denote by $\mathrm{wt}(c)$ the \textbf{Hamming weight} of $\mathbf{c}$, that is, the cardinality of $\mathrm{supp}(\mathbf{c})=\{i\mid c_i\neq 0\}$, the \textbf{support} of $\mathbf{c}$.
 The \textbf{minimum distance} of a code $\mathcal K\in R^n$, i.e. the minimum Hamming weight of the nonzero elements in $\mathcal K$, will be denoted by $d(\mathcal K)$.

\begin{definition}\label{def:socle} Let $R$ be a local ring with
maximal ideal $M=\hbox{rad}(R)$ and residue field $\mathbb F_q=\overline R$. The \textbf{socle}
$\mathfrak{S}(\mathcal{K})$ of an $R$-linear code $\mathcal K$ is defined as the sum of all its irreducible $R$-submodules.
\end{definition}

{Accordingly} to \cite{Nechaev99} the equality
\[
\mathfrak{S}(\mathcal{K})=\{\mathbf c\in\mathcal K\ |\ M\mathbf c=0\}
\]
holds for any $R$-linear code $\mathcal K$. So we may consider
$\mathfrak{S}(\mathcal{K})$  as a linear space over the field
$\mathbb F_q$ where $\bar r\cdot \mathbf{c}
=r\mathbf{c}$ for all $\bar r\in \mathbb F_q,\, \mathbf{c}\in
\mathfrak{S}(\mathcal{K})$.

\begin{lemma}\label{lem:Nechaev99}
Let $R$ be a local ring with maximal ideal $M$ and $\mathcal K$ an
$R$-linear code of length $n$. Then $\mathfrak{S}(\mathcal{K})$ is a
linear code of length $n$ over the field ${\mathbb F}_q=R/M$ and
$d(\mathcal K)=d(\mathfrak{S}(\mathcal{K}))$.
\end{lemma}
\begin{proof}
It is a direct translation of Proposition 5 in \cite{Nechaev99}.
\end{proof}

\begin{proposition}
    In the conditions of Theorem \ref{H} $d(\mathcal K)=d(\mathcal L)$, where $\mathcal L$ is the code $\left\langle \overline{G_1},\dots,\overline{G_t}\right\rangle +\overline I$ in $\mathbb F_q[X_1,\dots,X_r]/\left\langle \overline t_1(X_1),\dots,\overline t_r(X_r)\right\rangle$.
\end{proposition}
\begin{proof}
   The socle of the code $\mathcal K$ is
   $\mathfrak{S}(\mathcal K)=\left\langle a^{t-1}G_1,a^{t-1}G_2,\dots,a^{t-1}G_t\right\rangle+I$,
   that can be seen as a linear code over $\mathbb F_q$. Consider the $\mathbb F_q$-vector space isomorphism
   $\phi:a^{t-1}R[X_1,\dots,X_r]/I\to \mathbb F_q[X_1,\dots,X_r]/\overline I$, given by
   $a^{t-1}g+I\to \overline g+\overline I$ to conclude the result.
\end{proof}

In the general situation we can not state that the minimum distance of a semisimple code $\mathcal K$ is equal to the minimum distance of the code $\overline{\mathcal K}$. The more we can say is that, if $\overline{\mathcal K}\not=0$, then $d(\mathcal K)\le d(\overline{\mathcal K})$. However, there is one subclass of multivariable semisimple codes for which the equality holds.

\begin{definition}
In the conditions of Theorem \ref{H}, the code $\mathcal K$ is called \textbf{Hensel lift of a multivariable semisimple code} if $\left\langle G_1+I\right\rangle \not =I$ and $\left\langle G_i+I\right\rangle=0$, for all $i=2,\dots,t$.
\end{definition}

 This notion generalizes the definition of a \emph{Hensel lift of a cyclic code} introduced in \cite{Norton}. For this class of codes we have the following result.

\begin{corollary}
    If $\mathcal K\not =0$ is a Hensel lift of a multivariable semisimple code, then
   $d(\mathcal K)=d(\overline{\mathcal K})$.
\end{corollary}
\begin{proof}
    As noticed above the inequality $d(\mathcal K)\le d(\overline{\mathcal K})$ holds. On the other hand, since $\mathcal K$ is a Hensel lift of a multivariable semisimple code, we have that $\mathcal L=\overline{\mathcal K}$ and the result follows from the previous proposition.
\end{proof}

This collorary generalizes Collorary 4.3 in \cite{Norton} for Hensel lift of cyclic codes. Moreover, all classical bounds on distances for semisimple codes over fields (BCH, Hartmann-Tzeng, Roos, \dots ) also apply to their Hensel lifts. Remark that these bounds can be stated in the multivariable abelian case due to Proposition 8 in \cite{Poli}[Chapitre 6], {that we remind in Proposition \ref{prop:bound} below}.

\begin{definition}
    A multivariable semisimple code $\mathcal K\lhd R[X_1,\dots,X_r]/I$ is called \textbf{abelian}, if $I=\left\langle x_1^{e_1}-1,\dots,X_r^{e_r}-1\right\rangle$, where $e_1,\dots,e_r\in\mathbb N$.
\end{definition}

Let $S=\bigsqcup_{i=1}^{l}\bigsqcup_{j=1}^{s_i}C(\mu^{(i,j)})$ be the set of defining roots of a semisimple abelian code in $\mathbb F_q[X_1,\ldots, X_r]/\overline{I}$, where $C(\mu^{(i,j)})\in \mathcal C$
such that $p_{\mu^{(i,j)},1}=p_{\mu^{(k,l)},1}$ iff $i=k$. Consider for each class $C(\mu^{(i,j)})$ the polynomial:
\begin{eqnarray*}
 \frac{\bar t_1(X_1)}{p_{\mu^{(i,j)},1}(X_1)} & &\left( \prod_{k=2}^{r} \frac{\bar t_k(X_k)}{p_{\mu^{(i,j)},k}(X_k)}\prod_{k=2}^r \pi_{\mu^{(i,j)},k}(X_2,\ldots ,X_r)\right) =\\ \frac{\bar t_1(X_1)}{p_{\mu^{(i,j)},1}(X_1)} & &\left( F_{ij}(X_2,\ldots ,X_r)\right) \end{eqnarray*}

Here $p_{\mu^{(i,j)},k}\ k=1,\dots,r$, and $\pi_{\mu^{(i,j)},k}$ $k=2,\ldots ,r$ are as in Definition \ref{def:pw}, and $F_{ij}\in \mathbb F_q[X_2,\ldots ,X_r]$ is uniquely determined by the class $C(\mu^{(i,j)})$. Let us consider the field $\mathbb F^{(i)}=\mathbb F_q(X_1)/p_{\mu^{(i,1)},1}(X_1)$, and the code $J_i$ generated by  $\sum_{j=1}^{s_i}F_{ij}$ in the algebra $\mathbb F^{(i)}[X_2,\ldots , X_r]/\left\langle \bar t_2, \ldots , \bar t_r\right\rangle $, $i=1,\dots,l$.

\begin{proposition}\label{prop:bound}{
With the notations above, the minimum weight of a semisimple code over a field $\mathbb F_q$ and of the corresponding Hensel lift over $R$ is at least $\mathrm{min}_{1\leq i\leq l}\{d_i\cdot \delta_i\}$ where $d_i$ is the minimum weight of the  code in $\mathbb F_q[X_1]/\bar t(X_1)$ generated by
$$ \frac{\bar t(X_1)}{p_{\mu^{(i,1)},1}(X_1)\cdot\ldots\cdot p_{\mu^{(l,1)},1}(X_1)}$$
and $\delta_i$ is the minimum weight of the code $J_i$.}
\end{proposition}
\begin{proof} { It is a straight forward generalization of Lemma 3 and
 Proposition 8 in \cite{Poli}[Chapitre 6].}
\end{proof}

\begin{remark}
Notice that, in view of this result, the computation of the minimum distance of a semisimple abelian code in $r$ variables is reduced to computations of minimum distances of semisimple abelian codes in less number of variables.
\end{remark}

\section{Dual codes of abelian semisimple codes}

 In this section we describe the dual codes of abelian multivariable
 semisimple codes.
Notice that any defining ideal $I$ of abelian codes must satisfy the
following property: $(e_i,p)=1$, for all $i=1,\dots,r$, since the
code is semisimple. On the other hand, any semisimple abelian code
can be seen also as a \textbf{group code}, i.e., as an ideal of a
certain group ring. Namely, the group ring
$RG=R(\mathsf{C}_{e_1}\times\dots \times \mathsf{C}_{e_r})$, where
$\mathsf{C}_s$ is the cyclic group of order $s$.

\begin{definition}
Let $R[X_1,\dots,X_r]/I$ be a semisimple abelian code with
$I=\left\langle x_1^{e_1}-1,\dots,X_r^{e_r}-1\right\rangle$, then we
define the ring automorphism $\tau$ of $R[X_1,\dots,X_r]/I$ given by
$\tau(f(X_1,\dots,X_r))=f(X_1^{-1},\dots,X_r^{-1})=f(X_1^{e_1-1},\dots,X_r^{e_r-1})$.
It is clear that this automorphism preserves the Hamming weights of
a words.
\end{definition}

\begin{theorem}\label{dual}
If $\mathcal K=\left\langle  G_1,a G_2,\dots,a^{t-1}
 G_t\right\rangle+I$ is a semisimple abelian code in the conditions of Theorem \ref{H}, then
its dual code is $$\mathcal K^\perp=\left\langle  \tau(G_0),a \tau(G_t),\dots,a^{t-1}
 \tau(G_2)\right\rangle+I,$$ where the polynomials $\tau(G_i)$, $i=0,2,3,\dots,t$ are also in the conditions of Theorem \ref{H}.
\end{theorem}

\begin{proof}
    Let us first prove that $\mathcal K^\perp=\tau(\mathrm{Ann}(\mathcal K))$. For all $F+I\in
    R[X_1,\dots,X_r]/I$ we have that $F+I\in \tau(\mathrm{Ann}(\mathcal K))$ if, and only if, for all $Q+I\in \mathcal K$:
$$I=Q\tau(F)+I$$
$$=\sum_{i_1,\dots,i_r}q_{i_1,\dots,i_r}X_1^{i_1}\dots X_r^{i_r}
\sum_{j_1,\dots,j_r}f_{j_1,\dots,j_r}X_1^{e_1-j_1}\dots X_r^{e_r-j_r}+I$$
$$=\sum_{k_1,\dots,k_r}\left(\sum_{i_1,\dots,i_r}q_{i_1,\dots,i_r}f_{i_1-k_1\ (\hbox{mod }e_1),\dots,i_r-k_r\ (\hbox{mod }e_r)}\right)X_1^{k_1}\dots X_r^{k_r}+I$$
$$=\sum_{k_1,\dots,k_r}(\mathbf{q}\cdot\mathbf{z_{k_1,\dots,k_r}})X_1^{k_1}\dots X_r^{k_r}+I,$$
where $\mathbf{q}$ and $\mathbf{z_{k_1,\dots,k_r}}$ denote, respectively, the vector of coefficients of $Q$ and $X_1^{k_1}\dots X_r^{k_r}F$, in a fixed ordering of the terms in $R[X_1,\dots,X_r]/I$. Hence, $F+I\in \tau(\mathrm{Ann}(\mathcal K))$ if, and only if, for all $Q+I\in \mathcal K$ and for all $0\le k_1<e_1,\dots,1\le k_r<e_r$, $\mathbf{q}\cdot\mathbf{z_{k_1,\dots,k_r}}=0$, i.e., $\mathbf{y_{k_1,\dots,k_r}}\cdot\mathbf{f}=0$, where $\mathbf{y_{k_1,\dots,k_r}}$ denotes the vector of coefficients of $X_1^{-k_1}\dots X_r^{-k_r}Q$, that is if, and only if, $F+I\in \mathcal K^\perp$.

Notice that the polynomials $\tau(G_i)$, $i=0,\dots,t$ are in the conditions of Theorem \ref{H}, and so it is enough to see that $a^iG_{t+1-i}+I\in\mathrm{Ann}(\mathcal K)$, $i=0,\dots,t-1$,  to conclude the result
(here we denote $G_{t+1}=G_0$). Let $i,j=0,\dots,t-1$, if $i+j\ge t$, then $(a^iG_{t+1-i}+I)(a^jG_{j+1}+I)=
a^{i+j}(G_{t+1-i}G_{j+1})+I=I$ and, if $i+j<t $, then $\left\langle G_{t+1-i}+I\right\rangle\not=\left\langle G_{j+1}+I\right\rangle$, and so $(a^iG_{t+1-i}+I)(a^jG_{j+1})=I$, from the decomposition of $\mathcal K$ in Theorem \ref{H}.
\end{proof}

\begin{corollary}
    In the conditions of the previous theorem:
$$|\mathcal K^\perp|=|\bar R|^{\sum_{i=0}^{t-1}iN_i}$$
where $N_i$ is the number of zeros $\mu\in H_1\times\dots H_r$
of $\bar G_i, i=0,\dots, t-1$,
and
     $\mathcal K^\perp=\left\langle \tau(G_0)+a\tau(G_t)+\dots+a^{t-1}\tau(G_2)+I\right\rangle$
\end{corollary}

\begin{proof}
    The result follows from \cite[Proposition 2.11]{CN-FCR} and the fact that the polynomials $\tau(G_i)$ are in the conditions of Theorem \ref{H}.
\end{proof}

\begin{remark}
In view of Theorem \ref{dual} all the remarks concerning the distance of a code observed in the previous section can be applied also to its dual. Of course, the results about the minimum distance of a code and the minimum distance of its dual involving the MacWilliams identity for codes over Quasi-Frobenius modules \cite{Nechaev99} apply also in our case. In sake of brevity we will not get into details, though.
\end{remark}

\section{Self-dual abelian semisimple codes}

In the previous section we have described explicitly the dual code of a given abelian semisimple code $\mathcal K$. We want now to study conditions on $\mathcal K$ to be self-dual. Notice first
that, if the nilpotency index $t$ of $a$ is even, then there always exists a self-dual code, namely $\left\langle a^{\frac{t}{2}}\right\rangle$, that it is called the \textbf{trivial self-dual code}.
On the other hand, remember that any abelian code is also a group code and so the problem of existence of self-dual semisimple abelian codes can be reduced to the existence of self-dual group codes in $RG$.
This problem
 has been solved for some classes of rings $R$. In this direction an interesting work is \cite{Willems} where the existence of self-dual codes is characterized when $R$ is a Galois Ring. The techniques of proof make use of Group Representation Theory and can be also used when $R$ is a finite commutative chain ring. Namely, the following result holds.

\begin{theorem}
    Let $R$ be a finite chain commutative ring of characteristic $p$ with $a\in R$ such that $\left\langle a\right\rangle=\mathrm{rad}(R)$ with nilpotency index $t$, and let $G$ be a finite group. Then $RG$ contains a self-dual group code (that is, and ideal $\mathcal K\lhd RG$ such that $\mathbf{x}\cdot\mathbf{y}=0$, for all $x,y\in \mathcal K$) if, and only if, $p$ is odd and $t$ even, or $p$ and $t|G|$ are even.
\end{theorem}

\begin{proof}
    The proof is exactly the same that in the case of $R$ being a Galois Ring ({see \cite{Willems}}). This is due to the following two facts: any finite commutative chain ring $R$ is a Frobenius ring \cite{Wood}, and for any finite group $G$ we have a filtration
$$0 \subsetneq a^{t-1} RG \subsetneq \dots \subsetneq a^1 RG \subsetneq RG.$$
\end{proof}

In view of this result we can only expect to find non-trivial self-dual codes in the semisimple abelian case if,
and only if, $p$ and $|G|$ are even, or $t$ is even. The first case is clearly imposible, since
$|G|=\prod_{i=1}^re_i$ even implies that there exists some $e_i$ even and the code is not semisimple (notice that $p=2$). So we have only to study the case when $t$ is an even number. As a consequence to Theorem \ref{dual} we have the following result.

\begin{corollary}
Let $\mathcal K=\left\langle  G_1,a G_2,\dots,a^{t-1}
 G_t\right\rangle+I$ be a semisimple abelian code in the conditions of Theorem \ref{H}, then $\mathcal K$ is self-dual if, and only if, $\left\langle G_i+I\right\rangle=\left\langle \tau(G_j)+I\right\rangle$ when $i+j\equiv 1\ (\hbox{mod}\ t+1)$.
\end{corollary}

\begin{proof}
    By Theorem \ref{dual} we have $\mathcal K^\perp=\left\langle  \tau(G_0),a \tau(G_t),\dots,a^{t-1}
 \tau(G_2)\right\rangle+I$. Therefore, if $\left\langle G_i+I\right\rangle=\left\langle \tau(G_j)+I\right\rangle$ where
  $i+j\equiv 1\ (\hbox{mod}\ t+1)$, then $\mathcal K=\mathcal K^\perp$, and the code is self-dual.
Conversely, if $\mathcal K=\mathcal K^\perp$, then $\left\langle  G_1,a G_2,\dots,a^{t-1}
 G_t\right\rangle+I=\left\langle  \tau(G_0),a \tau(G_t),\dots,a^{t-1}
 \tau(G_2)\right\rangle+I$, and the result follows from the uniqueness of the ideals in Theorem \ref{H}.
\end{proof}

\begin{theorem}
    If $t$ is an even number, then there exist non-trivial self-dual semisimple abelian codes if, and only if, there
    exists $\mu\in H_1\times\dots\times H_r$ such that $C(\mu)\not=C(\mu^{-1})$, where
    $\mu^{-1}=(\mu_1^{-1},\dots,\mu_r^{-1})$.
\end{theorem}

\begin{proof}
    Let us first assume that there exists $\mu\in H_1\times\dots\times H_r$ such that $C(\mu)\not=C(\mu^{-1})$. Let $G+I$ be a generator of the semisimple abelian code
$\bigoplus_{\eta\not=\mu,\mu^{-1}}\left\langle h_\eta+I\right\rangle$ and consider:
$$\mathcal K=\left\langle a^{\frac{t}{2}-1}h_\mu,a^{\frac{t}{2}}G,a^{\frac{t}{2}+1}h_{\mu^{-1}}\right\rangle+I.$$
Since $\left\langle\tau(h_{\mu^{-1}})+I\right\rangle=\left\langle h_\mu+I\right\rangle$ and $\left\langle \tau(G)+I\right\rangle=\left\langle G+I\right\rangle$ we have, from the previous corollary, that
$\mathcal K$ is a non-trivial self-dual semisimple abelian code.

Conversely, if
$\mathcal K=\left\langle  G_1,a G_2,\dots,a^{t-1}
 G_t\right\rangle+I$ is a self-dual semisimple code, then for all $i,j$ such that $i+j\equiv 1\ (\hbox{mod }t+1)$ we have that $\left\langle G_i+I\right\rangle=\left\langle \tau(G_j)+I\right\rangle$. Assume now that $C(\mu)=C(\mu^{-1})$, for any $\mu\in H_1\times\dots\times H_r$. Then $\left\langle h_\mu+I\right\rangle =\left\langle h_{\mu^{-1}}+I\right\rangle=\left\langle \tau(h_\mu)+I\right\rangle$, and so $\left\langle G_j+I\right\rangle =
\left\langle \tau(G_j)+I\right\rangle=\left\langle G_i+I\right\rangle$, for all $i,j$ such that $i+j\equiv 1\ (\hbox{mod }t+1)$.
From the decomposition of Theorem \ref{H} we obtain that $\mathcal K=\left\langle  a^{\frac{t}{2}}+I\right\rangle$ is the trivial self-dual code.
\end{proof}

The existence of non-trivial self-dual codes can be eventually reduced to a number theoretical problem, as the following result shows.

\begin{corollary}
    If $t$ is an even number, then there exist non-trivial self-dual semisimple abelian codes if, and only if,
    $q^i\not\equiv-1\ (\emph{mod lcm}(e_1,\dots,e_r))$,
    for all na\-tu\-ral number $i$.
\end{corollary}

\begin{proof}
    From the previous theorem we have that non-trivial self-dual codes semisimple abelian codes do not exist if, and
    only if, $C(\mu)=C(\mu^{-1})$, for all $\mu\in H_1\times\dots\times H_r$. If $\xi_i$ denotes an
    $e_i$-th
    primitive root of unity, then this is equivalent to the condition for all $0\le a_i<e_i$, $i=1,\dots,r$, there exists a natural number $h$ such that $\xi_i^{-a_i}=\xi_i^{q^ha_i}$, i.e., $q^ha_i\equiv -a_i\ (\hbox{mod }(e_i))$. Therefore non-trivial self-dual codes do not exist if, and only if, there exists a natural number $h$ such that $q^h\equiv -1\ (\hbox{mod }(e_i))$ for all $i=1,\dots,r$, that is, $q^i\equiv-1\ (\hbox{mod lcm}(e_1,\dots,e_r))$.
\end{proof}

This result generalizes 4.4 Theorem in \cite{CN-FCR} for the case of self-dual cyclic codes. In this work it is also included a discussion about pairs of natural numbers $(q,n)$ for which $q^i\not\equiv-1\ (\hbox{mod }n)$,
    for all na\-tu\-ral numbers $i$, when $q$ is a prime number. The search of conditions for a pair of numbers to satisfy this property when $q$ is a power of a prime number is an open problem.

\end{document}